\documentclass[11pt]{amsart}

\usepackage{graphicx}
\usepackage{amscd}
\usepackage{amsmath}
\usepackage{amsfonts}
\usepackage{amssymb}
\usepackage{amsxtra}
\usepackage{xspace}
\usepackage{array}

\theoremstyle{plain}
\newtheorem{theorem}{Theorem}
\newtheorem{lemma}{Lemma}
\newtheorem{corollary}{Corollary}
\newtheorem{proposition}{Proposition}

\newtheorem{claim}{Claim}

\theoremstyle{definition}
\newtheorem{definition}{Definition}

\numberwithin{equation}{section}

\newcommand{\be}{\begin{enumerate}}
\newcommand{\ee}{\end{enumerate}}
\newcommand{\beq}{\begin{equation}}
\newcommand{\eeq}{\end{equation}}
\newcommand{\bprop}{\begin{proposition}}
\newcommand{\eprop}{\end{proposition}}

\newcommand{\complex}{\mathbb{C}}

\newcommand{\nat}{\mathbb{Z}_{\geq 0}}

\newcommand{\integers}{\mathbb{Z}}

\newcommand{\pfbegin}{\noindent {\em Proof:} }

\DeclareMathOperator{\len}{\ell}
\DeclareMathOperator{\supp}{supp} 
\DeclareMathOperator{\ch}{ch} 
\DeclareMathOperator{\htt}{ht} 
\DeclareMathOperator{\ct}{ct} 
\DeclareMathOperator{\Max}{Max}

\newcommand{\alm}[1][t]{a^\lambda_\mu(#1)} %no q
\newcommand{\lus}[1][t]{m^{\lambda}_{\mu}(#1)} 
\newcommand{\tko}[1][]{\mathcal{P}(#1; \,t)}
\newcommand{\kost}[1][\mu]{K_{\lambda #1}(t)}
\newcommand{\hl}[1][\lambda]{P_{#1}(t)}
\newcommand{\sch}[1][\lambda]{\chi_{#1}}
\newcommand{\scA}{\mathcal{A}}
\newcommand{\Ep}{\mathcal{E}}
\newcommand{\oh}{\mathcal{O}}

\newcommand{\kma}{\mathfrak{g}}
\newcommand{\csa}{\mathfrak{h}}
\newcommand{\bas}{L(\Lambda_0)}

\newcommand{\dtl}{\boldsymbol{\tilde{\varDelta}}}
\newcommand{\dtli}{\boldsymbol{\tilde{\varDelta}}^{im}}

\newcommand{\aoo}{\ensuremath{a^0_0(t)}\xspace}

\newcommand{\alolo}[1][t]{\ensuremath{a^{\Lambda_0}_{\Lambda_0}(#1)}\xspace}
\newcommand{\kmaf}{\overset{\circ}{\mathfrak{g}}}

\newcommand{\chke}{\hat{\boldsymbol{\mu}}}
\newcommand{\Kf}{\overset{\circ}{K}}
\newcommand{\Qf}{\overset{\circ}{Q}}
\newcommand{\Pf}{\overset{\circ}{P}}
\newcommand{\Lf}{\overset{\circ}{L}}
\newcommand{\rhof}{\overset{\circ}{\rho}}

\newcommand{\lf}{\bar{\lambda}}
\newcommand{\wf}{\overset{\circ}{W}}
\newcommand{\alf}{\bar{\alpha}}
\newcommand{\schf}[1][\lambda]{{\overset{\circ}{\chi}}_{#1}}
\newcommand{\delp}{{\overset{\circ}{\Delta}}_+}

\newcommand{\tsf}{$t$-string function\xspace}
\newcommand{\tsfs}{$t$-string functions\xspace}

\newcommand{\kos}[1][\mu]{\ensuremath{K_{\lambda,#1}(t)}\xspace}
\newcommand{\Ept}{\Ep_t}
\newcommand{\ro}{\Delta}
\newcommand{\rp}{\Delta_+}

\newcommand{\rmin}{\Delta_-}
\newcommand{\rre}{\Delta^{re}}
\newcommand{\rrep}{\Delta^{re}_+}
\newcommand{\rim}{\Delta^{im}}
\newcommand{\rimp}{\Delta^{im}_{+}}

\newcommand{\sa}[1][A]{|#1|}

\newcommand{\skma}{symmetrizable Kac-Moody algebra\xspace}
\newcommand{\skmas}{symmetrizable Kac-Moody algebras\xspace}

\newcommand{\dsum}{\displaystyle\sum}
\newcommand{\dprod}{\displaystyle\prod}

\newcommand{\highroot}{\theta} % Bourbaki, jrs notation

\begin{document}

\title[]{Kostka-Foulkes polynomials for symmetrizable Kac-Moody algebras}
\author{Sankaran Viswanath}
\address{Department of Mathematics\\
Indian Institute of Science\\
Bangalore 560012, India.}
\email{svis@math.iisc.ernet.in}

\begin{abstract}
We introduce a generalization of the classical Hall-Littlewood and
Kostka-Foulkes polynomials to all \skmas. We prove that these Kostka-Foulkes
polynomials coincide with the natural generalization of Lusztig's
$t$-analog of weight multiplicities, thereby extending a theorem of
Kato. For $\kma$ an affine Kac-Moody algebra, we define $t$-analogs of
string functions and use Cherednik's constant term identities to
derive explicit product expressions for them.
\end{abstract}

\maketitle

\section{Introduction}
The theory of Hall-Littlewood polynomials and Kostka-Foulkes
polynomials associated to a finite dimensional simple Lie algebra
$\kma$ is a classical subject with numerous connections to
representation theory, combinatorics, geometry and mathematical
physics. We refer to \cite{ramnelsen} and \cite{jrs} for nice surveys. 
Let $W$ be the Weyl group of $\kma$, $P$ its weight lattice and $P^+$
the set of dominant weights.  The
Hall-Littlewood polynomials $\hl$ and the Weyl characters $\sch[\mu]$
($\lambda, \mu \in P^+$) each form a $\complex[t]$-basis of the ring $\complex[t][P]^W$. The
Kostka-Foulkes polynomials $\kost$ are the entries of the transition
matrix between these two bases.

In this article, we introduce a generalization of the Hall-Littlewood
polynomials to all {\em \skmas} $\kma$. This is defined by
\beq\label{objint}
\hl := \frac{1}{W_\lambda(t)} \frac{\sum_{w \in W} (-1)^{\len(w)} w\left(e^{\lambda+\rho}
\prod_{\alpha \in \rp} \left( 1-t e^{-\alpha} \right)^{m_\alpha}\right)}{e^{\rho} \prod_{\alpha \in \rp} 
(1 - e^{-\alpha})^{m_\alpha}},
\eeq
where $W$ is the Weyl group and $\rp$ the set of positive roots of $\kma$
(see \S\ref{notnsec} for complete notation). This is a straightforward
generalization of the classical definition, taking into account both
real and imaginary roots (with multiplicities) of $\kma$. However, it
is not immediately clear that \eqref{objint} is well-defined since $W$
and $\rp$ are infinite sets in general. Our first task (proposition
\ref{jfl}) is to establish the well-definedness of $\hl$.  
The coefficient of a typical term $e^\mu$
on the right hand side of \eqref{objint} is a priori only a power
series in $t$, but it will in fact be shown to  be  a polynomial.

Many of the properties of the classical Hall-Littlewood polynomials
remain true for the general $\hl$. For instance $P_\lambda(0) = \sch$,
the character of the irreducible representation $L(\lambda)$ of
$\kma$. One can then define the Kostka-Foulkes polynomials $
K_{\lambda\mu}(t)$ for $\kma$ in the usual manner:
$$ \sch = \sum_{\substack{\mu \in P^+\\\mu\leq\lambda}} K_{\lambda\mu}(t) P_{\mu}(t).$$
Here, the sum on the right is typically infinite, 
$K_{\lambda\mu}(t) \in \integers[[t]]$ and $K_{\lambda\lambda}(t) =1$.

One of the fundamental facts concerning classical Kostka-Foulkes polynomials
is that they coincide with Lusztig's $t$-analog of weight
multiplicities \cite{lusztig} given by
$$ \lus := \sum_{w \in W} (-1)^{\len(w)} \tko[w(\lambda+\rho) -  (\mu+\rho)],$$ 
where $\tko$ is the $t$-analog of Kostant's partition function. In
this classical case, the
equality $\kost=\lus$ was proved by Kato \cite{kato}, and later by
R. Brylinski \cite{rkg} (see also \cite{jrs}). We prove this equality
(theorem \ref{mainthm}) in our more general setting of symmetrizable
Kac-Moody algebras. Here, $\lus$ for $\kma$ is obtained in the natural manner by
replacing $\tko$ with the $t$-analog of the (generalized) Kostant
partition function for $\kma$ (equation \eqref{introrefs}). Our proof
is based on an adaptation of ideas of Macdonald in \cite{igmgod}. For
the case of finite dimensional $\kma$, this provides an elementary proof
different from those of Kato and Brylinski. For general $\kma$, as 
nice consequences of this theorem, one obtains that $\kost \in
\integers[t]$ and $K_{\lambda\mu}(1) = \dim(L(\lambda)_\mu)$.

Next, we  specialize to the case that $\kma$ is an untwisted affine
Kac-Moody algebra. Given dominant weights $\lambda, \mu$, we consider
the generating function $\alm:= \sum_{k \geq 0} K_{\lambda,
  \mu-k\delta}(t) \; q^k         $. We call these {\em $t$-string
  functions} of $\kma$. When $t=1$, they reduce (up to a power of $q$)
to usual string functions. We  exhibit a surprising
connection of certain $\alm$ to Cherednik's famous {\em constant term
  identities}. In fact, we show that Cherednik's constant term
identities  of Macdonald and Macdonald-Mehta type \cite{dahamc,dmmc}
are the precise ingredients necessary to determine the \tsfs
associated to the trivial and basic representations of $\kma$.

Finally, we study a relation between Kostka-Foulkes polynomials for
affine Lie algebras and their classical counterparts, thereby
establishing some partial results concerning {\em positivity}.

The paper is organized as follows: \S\S
\ref{notnsec}-\ref{moredef} are concerned with the definitions, 
and \S \ref{mainsec} proves the equality $\kos =\lus$ in our setting.
In \S \ref{affreallystart}, \tsfs of levels 0,1 for affine Kac-Moody
algebras $\kma$ are computed. \S \ref{positi} and \S \ref{misc}
prove the partial positivity results and some assorted facts
concerning the $\hl$.

\vspace{1.5mm}
\noindent
{\em Acknowledgements:}  
 The motivation for this work arose from the
AIM workshop  on ``Generalized Kostka polynomials'' in July 2005. The author would 
like to thank the American Institute of Mathematics for its
hospitality and the workshop organizers (Anne Schilling and Monica
Vazirani) for their invitation to attend.
The author would like to especially thank John Stembridge for his inspiring 
survey talk \cite{jrs} and subsequent discussions. 

\vspace{1.5mm}
\noindent
{\em Postscript:} After this paper was completed, we were informed of
an earlier unpublished manuscript \cite{groj} of Ian Grojnowski,  which
overlaps with our present work. Together with \cite{fgt}, this
manuscript forms part of Grojnowski's larger program of extending
the notion of {\em geometric Satake isomorphism} to double loop groups.

Our approach to the combinatorial results in the present work is different from
and more complete than that of \cite{groj}. While our proof of the main theorem (theorem \ref{mainthm}) is based
 on ideas of Macdonald  \cite{igmgod}, the proof sketched in \cite{groj} attempts to directly generalize
  Brylinski's arguments for the finite dimensional case \cite{rkg}. This latter approach however presents
  ``well-definedness'' difficulties for arbitrary symmetrizable Kac-Moody algebras;  these are not fully resolved in \cite{groj}.

\section{Definitions}\label{notnsec}
We will use the notations of Kac's book \cite{kac}.
Let $\kma$ be a \skma with (finite dimensional) Cartan subalgebra
$\csa$. Let $\alpha_i, i=1 \cdots n$ be the simple roots of $\kma$ and
 $P$, $Q$, $P^+$, $Q^+$ be the weight lattice, the
root lattice and the sets of dominant weights and non-negative integer
linear combinations of simple roots respectively. We will denote the
Weyl group of $\kma$ by $W$ and let $(\cdot,\cdot)$ be a nondegenerate, $W$-invariant symmetric
bilinear form on $\csa^*$. Let $\ro$
(resp. $\rp$, $\rmin$) be the set of roots (resp. positive, negative
roots) of $\kma$. Given $\lambda, \mu \in \csa^*$, we say $\lambda \geq
\mu$ if $\lambda - \mu \in Q^+$. Given $\lambda \in \csa^*$, define
$D(\lambda) :=\{ \gamma \in \csa^*| \gamma \leq \lambda\}$.

Our goal here is to generalize the classical definition of  Hall-Littlewood polynomials $\hl$ 
to all \skmas. However, when passing to arbitrary Kac-Moody
algebras, we will have to deal with well-definedness issues concerning
infinite sums, products etc.

First, define $\Ept$ to be
the set of all series of the form
\begin{equation}\label{ser}
\sum_{\lambda \in \csa^*} c_\lambda(t) \, e^{\lambda},
\end{equation}
where each $c_\lambda(t) \in \complex[[t]]$ and $c_\lambda =0$ outside the
union of a finite number of sets of the form $D(\mu), \, \mu \in \csa^*$.
The $e^\lambda$ are formal exponentials, which obey the usual rules :
$e^0=1, e^{\lambda + \mu} = e^\lambda e^\mu$. Extending this
multiplication $\complex[[t]]$-bilinearly to all of $\Ept$ makes it into a commutative,
associative $\complex[[t]]$-algebra. We also let $\Ep$ denote the $\complex$-subalgebra of $\Ept$ 
consisting of series where all $c_\lambda \in \complex$.
Formal characters of $\kma$ modules
from category $\oh$ are elements of $\Ep$ \cite{kac}.

\subsection{}
Let $\rre$ (resp. $\rre_\pm$) denote the set of real roots of $\kma$
(resp. positive/negative real roots), and similarly $\rim$
(resp. $\rim_\pm$) be the corresponding subsets of imaginary
roots. For each $\alpha \in \ro$, let $m_\alpha$ denote the root
multiplicity of $\alpha$. Let $C
\subset \csa^*$ be the fundamental Weyl chamber and $X = \bigcup_{w \in W} wC$
 be the Tits cone. Let $\rho$ be a Weyl vector of $\kma$ defined by
 $(\rho, \alpha^\vee_i) =1$ for all $i=1\cdots n$, where as usual
 $\alpha^\vee_i = \frac{2\alpha_i}{(\alpha_i,\alpha_i)}$. Finally,
 given $\lambda \in P^+$, let $L(\lambda)$ be the integrable
 $\kma$-module with highest weight $\lambda$.

Given $\lambda \in P^+$, let $$f_\lambda := e^{\lambda+\rho}
\prod_{\alpha \in \rp} \left( 1-t e^{-\alpha} \right)^{m_\alpha} \; \in \Ept.$$
We can write  $f_\lambda = \sum_{\mu \in P} b_{\lambda \mu}(t) \,
e^\mu$; let $\supp (f_\lambda)$ denote the set of all $\mu$ for which
$b_{\lambda\mu}(t) \neq 0$. Each $\mu \in \supp(f_\lambda)$ can
be written as $\mu = \lambda +\rho - \sa$. Here $\sa :=\sum_{\alpha
  \in A} \alpha$ and
\begin{equation}
\label{adesc}
\begin{split}
&A \text{ is a finite multiset
of positive roots such that each } \alpha \in \rp  \\
& \text{ occurs at most } m_\alpha \text{ times in } A. 
\text{ For convenience, we will let } \scA \\
&\text{ denote the set of all such
  multisets}.
\end{split}
\end{equation}
A straightforward calculation gives $ b_{\lambda \mu}(t) = \sum (-t)^{\# A}$ where the
sum is over all multisets $A$ as described in \eqref{adesc} such that
$\lambda + \rho -\sa = \mu$.

Our first observation (the general version of the argument of \cite{rkg})
is :
\begin{lemma}\label{basic}
Let $\mu \in \supp(f_\lambda)$. Then 
\be
\item \label{fir} $\mu \in X$ (the Tits
cone)
\item \label{sec} For all $w \in W$, we have $w \mu \leq \lambda + \rho$.
\ee
\end{lemma}

\pfbegin
We only prove the second assertion, since it clearly implies the first.
Let $\mu = \lambda +\rho - \sa$ for some $A \in \scA$.
Let $L(\rho)$ denote the integrable highest weight representation of
$\kma$ with highest weight $\rho$. It is well-known \cite[Ex
  10.1]{kac} that the formal character $\ch L(\rho) := \sum_{\mu \in P}
\dim(L(\rho)_{\mu}) e^\mu$ of the module $L(\rho)$ 
is given by $e^\rho \prod_{\alpha \in \rp} (1 +
e^{-\alpha})^{m_\alpha}$.  Since $A \in \scA$, the element $\gamma :=
\rho -\sa$ is a weight of $L(\rho)$. Consequently, for all $w \in W,
w \gamma \leq \rho$. Now, $\lambda \in P^+$ implies that   $w\lambda \leq \lambda$ as
well. Since $\mu = \lambda + \gamma$, our assertion is proved. \qed

\subsection{} \label{hldef}
As in the classical finite dimensional situation, we now consider the
operator
$$ J := \sum_{w \in W} (-1)^{\len(w)} w,$$
where $\len(\cdot)$ is the length function on $W$. Given $f =
\sum_{\mu \in P} c_\mu (t) e^\mu \in \Ept$, we let $J(f) := \sum_{w
  \in W} \sum_{\mu \in P} c_\mu (t) (-1)^{\len(w)} e^{w\mu}$ whenever
the right hand side is defined. We make a couple of remarks about this
definition: 
\be
\item Let $\mu \in X$; then $J(e^\mu) := \sum_{w \in W} (-1)^{\len(w)}
  e^{w\mu}$ is clearly a well-defined element of $\Ept$ provided the
  stabilizer $W_\mu$ of $\mu$ in $W$ is finite. Further, $1 < \# W_\mu
  < \infty$ implies that  $J(e^\mu)=0$.
\item For arbitrary $f \in \Ep_t$ however, $J(f)$ may no longer be
  well-defined. As an example, let $f :=J(e^\rho) \in \Ep_t$. Then it
  is easy to see that  $J(f)$ is not defined.
\ee
Our main well-definedness result is:
\begin{proposition}\label{jfl}
For each $\lambda \in P^+$, $J(f_\lambda)$ is a well-defined element
of $\Ept$.
\end{proposition}
This proposition allows us to define our main objects of interest.
\begin{definition} \label{hlpolydef}
\textrm{
Given $\lambda \in P^+$, let $W_\lambda(t) : = \sum_{\sigma \in
  W_\lambda}  t^{\len(\sigma)} \; \in \complex[[t]]$. Then
$$\hl := \frac{1}{W_\lambda(t)}
\frac{J(f_\lambda)}{J(e^\rho)}$$ will be called
  the {\em Hall-Littlewood function} associated to $\lambda$.}
\end{definition}
\noindent
Observe, by the Weyl-Kac character formula that  
$J(e^\rho)^{-1} = e^{-\rho} \prod_{\alpha \in \rp} 
(1 - e^{-\alpha})^{-m_\alpha}$.  Here we
  interpret $(1 - e^{-\alpha})^{- m_\alpha}$ as the power series $(1
  + e^{-\alpha} + e^{-2\alpha} + \cdots)^{m_\alpha}$. Further since
  $W_\lambda(0)=1$, the inverse of $W_\lambda(t)$ is an element of 
  $\complex[[t]]$. It is thus clear from
  proposition \ref{jfl} that  $\hl \in \Ept$.

The rest of this subsection will be devoted to the proof of
proposition \ref{jfl}.
Recall that $f_\lambda = \sum_{\mu \in P} b_{\lambda \mu}(t) \,
e^\mu$. As a first step toward showing well-definedness of
$J(f_\lambda)$, we prove :
\begin{claim}\label{clm}
$\mu \in \supp (f_\lambda)$ implies that  $\# W_\mu <\infty$.
\end{claim}
{\em Proof:} Write $\mu = \lambda + \rho - \sa$ with $A$ a fixed
multiset in $\scA$.  Then $w \in W_\mu$
implies that $w(\lambda + \rho - \sa) =  \lambda + \rho - \sa$ or
equivalently that $$ \left[ \rho - w\rho +   w(\sa) \right ] +
(\lambda - w \lambda) = \sa$$ For an element $\beta = \sum_{i=1}^n a_i
\alpha_i \; \in Q$, the ``height'' of beta is defined to be
$\htt(\beta):=\sum_i a_i$. Since $\lambda \geq w \lambda$, we have
$\htt( \rho - w\rho +  w(\sa)) \leq \htt(\sa)$. Lemma \ref{future}
(below) implies that $| \len(w) - \# A| \leq \htt(\sa)$ and hence that
$\len(w) \leq \htt(\sa) + \#A$. Since $A$ was fixed to begin with, this
implies that  $\ell(w)$ is bounded above. Thus there can only be finitely
many elements in $W_\mu$. \qed

We now state and prove the lemma that was referred to in the above
proof. This simple lemma will come into play again and enable us to
complete the proof of proposition \ref{jfl}.

\begin{lemma}\label{future} Let $w \in W$ and $A \in \scA$.
 Then $$ \htt(\rho - w \rho + w(|A|)) \geq |\len(w)
  - \# A|.$$
\end{lemma}
\pfbegin
Let $S(w) :=\{\beta \in \rp: w^{-1} \beta \in \rmin\}$ and put $A_{re}
:=\{ \alpha \in A: \alpha \in \rre\}$ and  $A_{im}
:=\{ \alpha \in A: \alpha \in \rim\}$. We recall that $\rho - w\rho =
\sum_{\beta \in S(w)} \beta$ and $\# S(w) = \len(w)$. Thus 
\begin{equation}
\label{rightabove}
\rho - w\rho + w(|A|) = \sum_{\beta \in S(w)} \beta +  \sum_{\substack{\alpha
    \in A_{re} \\ w\alpha \in \rmin}} w \alpha + \sum_{\substack{\alpha
    \in A_{re} \\ w\alpha \in \rp}} w \alpha + \sum_{\alpha \in A_{im}}
w \alpha.
\end{equation}
We observe that (i) $\alpha \in A_{im}$ implies that  $w\alpha \in \rp$
(ii) $\alpha \in A_{re}$ such that $w\alpha \in \rmin$ implies that $-w\alpha
\in S(w)$. Thus the terms in the second sum on the right hand side of
\eqref{rightabove} cancel with their negatives in the first sum. Let
the number of terms in the second sum be denoted $k(w,A)$. Clearly
$$k(w,A) \leq \min \{\len(w),\,\# A_{re}\} \leq \min \{\len(w),\,\# A\} =
\frac{1}{2} (\len(w) + \,\# A -  |\len(w) - \,\# A| ).$$ Thus, after this
cancellation, the right hand side of \eqref{rightabove} is a sum of
$r$ positive roots, where $r = \len(w) + \# A  - 2k(w,A) \geq
|\len(w) - \# A|$. Since each positive root has height at least 1,
our lemma is proved. \qed

We now complete the proof of proposition \ref{jfl}. Since 
$f_\lambda = \sum_{\mu \in P} b_{\lambda \mu}(t) e^{\mu}$,
we have $J(f_\lambda) = \sum_{\mu \in P} b_{\lambda \mu}(t)
J(e^{\mu})$. Note that each $J(e^{\mu})$ is well-defined by claim
\ref{clm}. To show the well-definedness of the entire sum, observe that 
$$ J(f_\lambda) = \sum_{w \in W} \sum_{A \in \scA} (-1)^{\len(w)}
(-t)^{\# A} e^{w(\lambda + \rho - |A|)}.$$
This is now a well-defined element of $\Ept$ if for each term
$e^\gamma$ that occurs on the right hand side, its coefficient is a
well-defined power series in $t$. In other words, for fixed $\gamma
\in P$ (by lemma \ref{basic}, we can assume $\gamma \leq \lambda +
\rho$) and $p \geq 0$, we need to show that the set $$\{(w,A) \in W
\times \scA: \# A =p \text{ and } w(\lambda + \rho - |A|) = \gamma
\}$$ is finite. Given $(w,A)$ in the above set, we have $\lambda +
\rho - \gamma = (\lambda - w \lambda) + (\rho - w\rho +
w(|A|))$. By lemma \ref{future}, $|\len(w) - \#A| \leq \htt(\lambda +
\rho - \gamma)$. So, $\len(w) \leq p + \htt(\lambda + \rho - \gamma)$,
which means that only finitely many $w$'s are possible. Fix one such $w$;
then $|A| = \lambda + \rho - w^{-1} \gamma$ is also fixed. Since the
elements of $A$ are positive roots, the number of possibilities for
$A$ corresponding to this $w$ is also finite. This completes the proof
of proposition \ref{jfl}. \qed

\section{Kostka-Foulkes polynomials}\label{moredef}
\subsection{}
We refer back to the definition of the Hall-Littlewood function $\hl$ in
\S\ref{hldef}. Observe that when $t=0$, $\hl$ reduces
to $J(e^{\lambda+\rho})/J(e^\rho) =: \sch$, the formal character of
the integrable highest weight module $L(\lambda)$. Given $\gamma \in
X$ such that $\#W_\gamma < \infty$, consider $J(e^\gamma)/J(e^\rho)$;
this equals $\epsilon_\gamma \sch[{[\gamma]-\rho}]$ where $[\gamma]$ is
the unique dominant weight in the Weyl group orbit of $\gamma$, and
$\epsilon_{\gamma} = 0$ if $\# W_\gamma > 1$ and equals
$(-1)^{\len(\sigma)}$ if $\# W_\gamma = 1$ and $\sigma \in W$ is such
that $\sigma \gamma = [\gamma]$. This discussion implies that
$\hl= W_\lambda(t)^{-1}  \sum_{\gamma \in P} b_{\lambda \gamma}(t)
\epsilon_\gamma \sch[{[\gamma]-\rho}]$. It is clear then that there exists 
$c_{\lambda \mu}(t) \in \complex[[t]]$ such that
\beq\label{plclchim}
 \hl = \sum_{\mu \in P^+} c_{\lambda\mu}(t) \sch[\mu]
\eeq
The right hand side is an infinite sum in general. We also get: 
\begin{align}
c_{\lambda\mu}(t) &=  W_{\lambda}(t)^{-1} \sum_{w \in W}
(-1)^{\len(w)} \, b_{\lambda,w^{-1}(\mu+\rho)}(t)\notag\\
&= W_{\lambda}(t)^{-1} \sum \sum (-1)^{\len(w)} (-t)^{\# A},  \label{corr9}
\end{align}
where the double sum in the last equation is over all $w \in W$ and
$A \in \scA$ such that $w(\lambda + \rho - |A|) = \mu + \rho$.

We have the following easy facts: (i) $c_{\lambda\mu}(t) \neq 0$
implies that  $\mu \leq \lambda$. This is an easy consequence of lemma
\ref{basic}. Also see \S \ref{misc} for a strengthening of this
fact. (ii) $c_{\lambda \lambda}(t)=1$. This calculation can be
performed by suitably modifying R. Brylinski's argument
in \cite[2.7]{rkg} for the finite case.

Thus, the change of coordinate matrix from the $\hl$ to the
$\sch[\mu]$ is an infinite rank,
 lower triangular matrix (after suitably ordering  the $\lambda
 \in P^+$) with ones on the diagonal. Inverting
the matrix, we can write
\beq\label{chikpl}
 \sch = \sum_{\mu \in P^+} K_{\lambda\mu}(t) P_{\mu}(t).
\eeq
By analogy to the classical case, we will call $K_{\lambda\mu}(t)$
{\em Kostka-Foulkes polynomials} for the \skma $\kma$. At the moment, we only
know that $K_{\lambda\mu}(t)$ is a power series in $t$, but theorem
\ref{mainthm} will establish that $\kost \in \integers[t]$. Observe
that $\kost \neq 0$ implies that  $\mu \leq \lambda$, and $\kost[\lambda]=1$.

\subsection{}
We now derive an alternative expression for $\hl$. If $w \in W$,
recall that $S(w) =\{\beta \in \rp: w^{-1}\beta \in \rmin\}$. Given $\alpha
\in \rp$, either (i) $w \alpha \in \rp$, in which case $\beta
:= w\alpha \not\in S(w)$ or (ii) $w \alpha \in \rmin$ in which case
$\beta:= -w \alpha \in S(w)$. We also have (iii) $w \rho - \rho = - \sum_{\beta
  \in S(w)} \beta$. So
\begin{align*}
w(f_\lambda) &= e^{w(\lambda+\rho)} \prod_{\alpha \in \rp} 
(1-te^{-w \alpha})^{m_\alpha} \\
& = (-1)^{\len(w)} e^{w\lambda} e^\rho \prod_{\substack{\beta \in \rp
    \\ \beta \not\in S(w)}}  (1-te^{-\beta})^{m_\beta} \prod_{\beta \in S(w)}
(t- e^{-\beta}).
\end{align*}
Here, we used the additional facts that (iv) $\len(w) = \# S(w)$, (v) $\beta \in
S(w)$ implies that  $\beta \in \rre$ and so $m_\beta =1$ and finally (vi)
$m_\alpha = m_{w\alpha}$ for all $w\in W$. Together with
the Weyl-Kac denominator formula $J(e^\rho) = e^\rho \prod_{\alpha
  \in \rp} (1 - e^{-\alpha})^{m_\alpha}$ and definition \ref{hlpolydef}, this gives
\begin{equation}\label{aliterpl}
\hl = \frac{ \displaystyle\sum_{w \in W} \left[ \,e^{w\lambda}
  \displaystyle\prod_{\substack{\beta \in \rp \\ \beta \not\in S(w)}}
  (1-te^{-\beta})^{m_\beta} \displaystyle\prod_{\beta \in S(w)} 
(t- e^{-\beta}) \,\right]}{W_\lambda(t) \; \displaystyle\prod_{\alpha \in \rp} (1-e^{-\alpha})^{m_\alpha}}.
\end{equation}

\section{The main theorem}\label{mainsec}
\subsection{}
Our main goal in this section will be to establish the relation
between the Kostka-Foulkes polynomial  $\kost$ for $\kma$ and Lusztig's
$t$-analog of weight multiplicities. We use an adaptation of ideas  of Macdonald
\cite{igmgod}. First define the element $\dtl \in \Ept$ by
\beq \label{vardeldefn}
\dtl := \frac{\prod_{\alpha \in \rp}
    (1-e^{-\alpha})^{m_\alpha}}{\prod_{\alpha \in \rp}
      (1-t e^{-\alpha})^{m_\alpha}} =  \prod_{\alpha \in \rp}
    \left[(1-e^{-\alpha})(1+t e^{-\alpha} + t^2 e^{-2\alpha} +
      \cdots)\right]^{m_\alpha}.
\eeq
Using equation \eqref{aliterpl}, we have
\begin{align} 
\dtl \hl &= \frac{\displaystyle\sum_{w \in W}
 e^{w\lambda} \displaystyle\prod_{\substack{\beta \in \rp \\ \beta \not\in S(w)}}
 (1-te^{-\beta})^{m_\beta} \displaystyle\prod_{\beta \in S(w)} (t-
 e^{-\beta})}{W_\lambda(t) \; \displaystyle\prod_{\alpha \in \rp} (1-t e^{-\alpha})^{m_\alpha}} \notag\\
&= \frac{1}{W_\lambda(t)} \sum_{w \in W} e^{w\lambda} \prod_{\alpha
      \in S(w)} \frac{t-e^{-\alpha}}{1-te^{-\alpha}}. \label{dlpl}
\end{align}
Clearly $\dtl \hl$ is also an element of $\Ept$.
Our main observation concerning \eqref{dlpl} is the following:
\begin{proposition} \label{dellm}
Let $\mu \in P^+$. The coefficient of $e^\mu$ in $\dtl
\hl$ equals $\delta_{\lambda\mu}$ (i.e is 1 if $\lambda = \mu$ and 0
otherwise).
\end{proposition}
\pfbegin
For convenience, we define (for $w \in W$):
$$ p_w := e^{w\lambda} \prod_{\alpha
      \in S(w)} \frac{t-e^{-\alpha}}{1-te^{-\alpha}} =  e^{w\lambda} \prod_{\alpha
      \in S(w)} (t-e^{-\alpha})(1 + te^{-\alpha} + t^2 e^{-2\alpha}+\cdots).$$
Thus $\dtl \hl = W_\lambda(t)^{-1} \sum_{w \in W} p_w$.
Now, suppose the term $e^\mu$ occurs in  $p_w$ for some $w \in W$. Then
$\mu$ can be written as $$\mu = w\lambda - \sum_{\alpha \in S(w)}
n_\alpha \alpha$$ with $n_\alpha \in \integers^{\geq 0}$. Thus $\mu
\leq w\lambda$. On the other hand, applying $w^{-1}$ we have:
$$w^{-1}\mu = \lambda - \sum_{\alpha  \in S(w)} n_\alpha
 w^{-1}(\alpha) \geq \lambda$$
since $\alpha \in S(w)$ implies that $w^{-1}(\alpha) \in \rmin$. Since $\mu \in P^+$
we also have $\mu \geq w^{-1}\mu$. Putting all these together, we get
$ \lambda \leq w^{-1}\mu \leq \mu \leq w \lambda$. Finally $\lambda \in
P^+$ implies that $\lambda \geq w \lambda$ as well, thereby forcing 
\begin{equation}\label{harvest}
  \lambda = w^{-1}\mu = \mu = w \lambda.
\end{equation}

\noindent
\underline{Consequences} (i) If $\mu \neq \lambda$ : the coefficient
of $e^{\mu}$ in $p_w$ is zero for all $w \in W$. (ii) If $\mu = \lambda$,
then by equation \eqref{harvest}, $w$ must satisfy $\lambda = w
\lambda$ i.e, $\lambda \in W_\lambda$. In this case, the coefficient
of $e^\lambda$ in $p_w$ is clearly $t^{\# S(w)} =
t^{\len(w)}$. Thus the coefficient of $e^\lambda$ in
$\dtl \hl$ is $W_\lambda(t)^{-1} \sum_{w \in W_\lambda}
t^{\len(w)} =1$. \qed

\vspace{2mm}
\noindent
{\bf Remark:} Setting $t=0$ in proposition \ref{dellm}, we have $\hl = \sch$, $\dtl =
\prod_{\alpha \in \rp}(1-e^{-\alpha})^{m_\alpha}$. In this case,  we
recover the well-known fact that the coefficient of
$e^{\mu + \rho}$ in $J(e^{\lambda + \rho})$ is $\delta_{\lambda\mu}$.

\subsection{} 
Let $\tko$ denote the  $t$-analog of the (generalized) Kostant
partition function for $\kma$. This is  defined via
the relation
\beq\label{introrefs}
\frac{1}{\prod_{\alpha \in \rp} (1-te^{-\alpha})^{m_\alpha}} =:
\sum_{\gamma \in Q^+} \tko[\gamma] \: e^{-\gamma}.
\eeq
Let $\lambda \in P^+$ and $\mu \in P$. Then Lusztig's $t$-analog of
weight multiplicity which we will denote by $\lus$, is defined to
be
$$\lus := \sum_{w \in W} (-1)^{\len(w)} \tko[w(\lambda+\rho) -
  (\mu+\rho)].$$ 
The right hand side has only finitely many nonzero terms and hence
  $\lus \in \integers[t]$. It is a theorem of Kato \cite{kato}
that for a finite dimensional simple Lie algebra $\kma$, we have
  $\kost = \lus$ for all $\lambda, \mu \in P^+$. Our main theorem
  extends this to all \skmas:
\begin{theorem}\label{mainthm}
Let $\kma$ be a \skma and $\lambda, \mu \in P^+$. Then $\kost = \lus$.
\end{theorem}
\pfbegin
We will employ the standard notation $[e^{\mu}] \; f$ to
denote the coefficient of $e^\mu$ in an expression $f$. Using equation
\eqref{chikpl}, we have
\begin{align*}
[e^{\mu}] \; \dtl \sch &= 
\sum_{\gamma \in P^+ } \kost[\gamma] \;
([e^{\mu}] \; \dtl \,\hl[\gamma]) \\
&= \kost,
\end{align*}
where for the last equality, we used proposition \ref{dellm}. Now, the
Weyl-Kac character formula gives
\begin{align}
\dtl \sch &= \frac{\sum_{w \in W} (-1)^{\len(w)}
  e^{w(\lambda+\rho) - \rho}}{\prod_{\alpha \in
  \rp}(1-te^{-\alpha})^{m_\alpha}} \notag\\
&= (\,\sum_{w \in W} (-1)^{\len(w)} e^{w(\lambda+\rho) - \rho}\,)
(\,\sum_{\gamma \in Q^+} \tko[\gamma] \: e^{-\gamma}\,). \label{lasst}
\end{align}
Using \eqref{lasst}, a direct calculation gives that the coefficient of
$e^\mu$ in $\dtl \sch$ is $\sum_{w \in W}
(-1)^{\len(w)} \tko[w(\lambda+\rho) -   (\mu+\rho)]$, which is
precisely $\lus$. \qed

\vspace{2mm}
\noindent
{\bf Remark:} For finite dimensional $\kma$, Ranee Brylinski
\cite{rkg} gave a proof of this theorem that was more elementary
than Kato's original proof. Ours is yet another  proof of this theorem,
which works just as well for all \skmas.

\begin{corollary} \label{tanalog}
Let $\kma$ be a \skma and $\lambda, \mu \in
P^+$. Then 
\be
\item $\kost \in \integers[t]$ and
\item $K_{\lambda \mu}(1)=\lus[1] = \dim(L(\lambda)_\mu)$.
\ee
\end{corollary}

\noindent
{\em Proof:} The first part follows from the fact that $\lus$
is an element of $\integers[t]$. Recall that a priori, we could only
say that $\kost$ was a power series in $t$.
The second part is a direct consequence of Kostant's weight multiplicity formula. \qed

On account of theorem \ref{mainthm}, we will often refer to the $\kos$
as $t$-weight multiplicities.

\section{Kostka-Foulkes polynomials for affine Kac-Moody algebras}
\label{affreallystart} 

We now specialize to the case that $\kma$ is an untwisted affine Kac-Moody
 algebra.  Let $\kmaf$ be the underlying finite dimensional simple Lie
 algebra of rank $l$, say. If $\delta$ is the null root, then
 $\rimp =\{k\delta: k \geq 1\}$ and each imaginary root has multiplicity
$l$. 

Suppose $\lambda, \mu$ are 
dominant weights of $\kma$ such that $\mu \leq \lambda$. We would like to
study $\kos$; observe here that for each $k \geq 0$, $\mu - k\delta$
is also a dominant weight $\leq \lambda$.
Let $\Max(\lambda):=\{\mu \in P^+: \mu \leq
 \lambda; \, \mu + \delta \nleq \lambda\}$
 (note: this notion is slightly different from that of Kac \cite[\S
 12.6]{kac} in that we do not require $\mu$ to be a weight of the
 $\kma$-module $L(\lambda)$). For each $\mu \in \Max(\lambda)$, we
 form the generating function 
\beq
\alm := \sum_{k \geq 0} K_{\lambda, \mu-k\delta}(t) \; e^{-k\delta}
\eeq 
of $t$-weight multiplicities along the $\delta$-string through
$\mu$. We will find it convenient to let $q:=e^{-\delta}$ in the rest
of the paper; thus for example, $\alm =  \sum_{k \geq 0} K_{\lambda,
  \mu-k\delta}(t) \; q^k$. When
$t=1$, $\alm[t]$ reduces to the generating function for ordinary weight
multiplicities along the $\delta$-string. Thus $\alm[1]$ 
is  (up to multiplication by a power of $q$) a {\em
 string function} for the module  $L(\lambda)$ \cite{kac}. 
By mild abuse of  terminology, we will call $\alm$ a {\em $t$-string
  function}.

We now recall the
definition of the {\em constant term} map $\ct(\cdot)$ from
\cite{igmgod}. Given $f = \sum_{\lambda}
f_{\lambda} e^\lambda \in \Ept$, define $\ct(f):=\sum_{k \in
  \integers} f_{k\delta} \: e^{k\delta}$. If $g \in \Ept$ is such that
$g=\ct(g)$, we recall the observation of \cite[(3.1)]{igmgod} that
$\ct(fg) = g \ct(f)$ for all $f \in \Ept$.

\subsection{\tsf of level 0}
The simplest \tsf is that of level 0. It is enough to consider 
$\lambda=0$; in this case $\Max(0) = \{0\}$.
The \tsf \aoo is easy to describe. Observe that theorem
 \ref{mainthm}  implies  that $K_{0,-k\delta}(t)$ is the
coefficient of $e^{-k\delta}$ in 
$$\frac{\sum_{w \in  W} (-1)^{\len(w)} e^{w\rho -\rho}}
{\prod_{\alpha \in \rp} (1-te^{-\alpha})^{m_\alpha}} = 
\frac{\prod_{\alpha\in \rp} (1-e^{-\alpha})^{m_\alpha}}
{\prod_{\alpha\in \rp} (1-te^{-\alpha})^{m_\alpha}} = \dtl$$
(notation from \eqref{vardeldefn}). In other words, we have
\beq\label{aooeq}
\aoo = \ct( \dtl).
\eeq

Next, we shall  show that equation \eqref{aooeq} is equivalent to an identity
due to Macdonald \cite[(3.8)]{igmgod}. Since $w\delta = \delta$
for all $w \in W$, equation \eqref{aliterpl} gives 
\beq \label{simpobs}
P_{\lambda + c \delta}(t) = e^{c\delta} P_\lambda(t) \;\; \forall \lambda \in P^+, c \in \complex.
\eeq
As a consequence, $1 = \sch[0] = \aoo \, \hl[0]$. From equation
\eqref{dlpl}, we have $W(t) \,\dtl \, \hl[0] =  \sum_{w \in W} \prod_{\alpha
      \in S(w)} \frac{t-e^{-\alpha}}{1-te^{-\alpha}}$. Let $\chke:=\prod_{\alpha \in \rrep}
\frac{1-e^{-\alpha}}{1-te^{-\alpha}}$ and $\dtli := \prod_{k \geq 1}
\left(\frac{1-e^{-k\delta}}{1-te^{-k\delta}}\right)^{l}$.
Observe that (i) $\ct(\dtli) = \dtli$ and hence, (ii)
$\ct(\dtl) = \ct(\chke \, \dtli) = \dtli \, \ct(\chke)$. 
Putting all this together and using \eqref{aooeq}, we deduce the following result.
\begin{proposition} \label{macd}
\beq\label{macid}
\frac{1}{W(t)} \,\sum_{w \in W} \prod_{\alpha
      \in S(w)} \frac{t-e^{-\alpha}}{1-te^{-\alpha}} = \frac{\chke}{\ct(\chke)}.
\eeq
\end{proposition}
This is essentially identity (3.8) of Macdonald \cite{igmgod} (in the
case where all $t_\alpha$ are equal). In fact, \cite{igmgod} was the
starting point for much of the present article. For many interesting
consequences of this identity, see \cite{igmgod}.

Finally, we recall the following (special case of an) 
identity due to Cherednik \cite{dahamc}:
\beq\label{cher1}
\ct(\chke ) =  \prod_{\alpha \in \delp} \frac{(t^{(\rhof,
    \alpha^{\vee})} \,q; \, q)_\infty^2}{(t^{(\rhof,    \alpha^\vee) +1} 
\,q; \, q)_\infty (t^{(\rhof,    \alpha^\vee) - 1} 
\,q; \, q)_\infty},
\eeq
where we have used the usual 
notation $(x;q)_\infty := \prod_{n=0}^{\infty} (1-xq^n)$. 
Here $\rhof$, $\delp$ are the Weyl vector and set of positive roots
(resp.) of $\kmaf$ and $\alpha^\vee = 2\alpha / (\alpha,\alpha)$. This
allows us to write \eqref{aooeq} and \eqref{macid} as explicit
infinite products (cf (3.14) of \cite{igmgod}). 

\subsection{The level 1 \tsf}
Let $\bas$ be the {\em basic
 representation} of $\kma$. Here $\Lambda_0$ is
the fundamental weight corresponding to the extended (zeroth) node of
the Dynkin diagram of $\kma$.  From \cite[Chap 12]{kac}, we  again
 have $\Max(\Lambda_0) = \{\Lambda_0\}$. We would  like to study
 the \tsf $\alolo$.

Now, suppose $\kma$ is one of the {\em simply-laced}  affine
Kac-Moody algebras $A_l^{(1)}, D_l^{(1)}, E_{6/7/8}^{(1)}$. In this case,
the  multiplicity of the weight $\Lambda_0 - k \delta$ ($k \geq 0$) in the
  basic representation $\bas$ equals $p_l(k)$, the number of
  partitions of $k$ into parts of $l$ colors
  \cite[prop. 12.13]{kac}. Since $K_{\Lambda_0, \Lambda_0 - k
    \delta}(1)$ is equal to this weight multiplicity,     this means
\beq \label{strid}
a_{\Lambda_0}^{\Lambda_0}(1) = \frac{1}{(q;q)_{\infty}^l}.
\eeq
In light of identity \eqref{strid}, we can now ask if
$a_{\Lambda_0}^{\Lambda_0}(t)$ has a nice closed form expression?
We have the following theorem:
\begin{theorem}\label{thm-lev1-tsf} 
Let $\kma$ be an untwisted affine Kac-Moody algebra. Then
\beq \label{thmeq}
\frac{a_{\Lambda_0}^{\Lambda_0}(t)}{a_{\Lambda_0}^{\Lambda_0}(1)}=
\dprod_{i=1}^l \frac{(q;\,q)_{\infty}}{ (t^{d_i} q;\,q)_{\infty}} =
   \dprod_{i=1}^l \dprod_{n=1}^\infty \frac{1-q^n}{1-t^{d_i} q^n}
\eeq
where $d_1, d_2, \cdots, d_l$ are the degrees (exponents plus one) 
of  the underlying finite  dimensional simple Lie algebra $\kmaf$.
\end{theorem}
\pfbegin
Let $\lambda \in P^+$; using \eqref{simpobs}, we can
 write $\sch = \sum_{\gamma
  \in \Max(\lambda)} a^{\lambda}_{\gamma}(t)\, \hl[\gamma]$. For 
$\mu \in \Max(\lambda)$, it is easy to see from proposition
\ref{dellm} that $\ct(e^{-\mu}\, \dtl\, \hl[\gamma]) =
\delta_{\gamma,\mu}$, and hence
\beq\label{almeq}
\alm = \ct(e^{-\mu}\, \dtl\, \sch).
\eeq
Applying this to our situation gives $a_{\Lambda_0}^{\Lambda_0}(t) =
\ct(e^{-\Lambda_0}\, \dtl\, \sch[\Lambda_0])$. We now let 
$$\Theta :=\sum_{\alpha \in \Qf} e^\alpha  q^{-(\alpha,\alpha)/2}$$
be the theta function of the root lattice $\Qf$ of $\kmaf$.
We then have \cite[lemma 12.7]{kac} 
\beq
e^{-\Lambda_0} \, \chi_{\Lambda_0} = a_{\Lambda_0}^{\Lambda_0}(1) \,
\Theta. 
\eeq
Recalling $\dtl = \chke \, \dtli$ and using the
aforementioned properties of the constant term map gives
\beq \label{puttog}
\frac{a_{\Lambda_0}^{\Lambda_0}(t)}{a_{\Lambda_0}^{\Lambda_0}(1)} =
\dtli \, \ct(\chke \, \Theta).
\eeq
We recognize the last term from Cherednik's work on the ``difference
analog of the Macdonald-Mehta conjecture'' \cite{dmmc}.
Specifically, adjusting for our sign convention and taking all $q_\alpha$ equal in \cite[(5.10)]{dmmc}, we have :
\begin{theorem} (Cherednik) \label{cher}
$$\ct(\chke \, \Theta) =  \prod_{\alpha \in \delp} \frac{(t^{(\rhof,
    \alpha^{\vee})} \,q; \, q)_\infty}{(t^{(\rhof,    \alpha^\vee) +1} 
\,q; \, q)_\infty}
$$
where $\rhof$, $\delp$ are the Weyl vector and set of positive roots
(resp.) of $\kmaf$ and $\alpha^\vee = 2\alpha / (\alpha,\alpha)$.
\end{theorem}

\noindent
For each $j\geq 1$, let $p_j$ denote the number of $\alpha \in \delp$
for which $(\rhof,\alpha^\vee)=j$ (i.e the number of {\em height} $j$
coroots). We have the following well-known fact \cite{kostant3d}:
$$p_j - p_{j+1} = \text{ number of times } j \text{ occurs as an
  exponent of } \kmaf. $$
Let $e_1, e_2, \cdots, e_l$ be the exponents of $\kmaf$. These are
related to the degrees by $d_i = e_i+1$.
Theorem \ref{cher} then gives $$\ct(\chke \, \Theta) = (tq; \,
  q)_\infty^l \; \prod_{i=1}^l \frac{1}{(t^{e_i+1}\,q; \,
  q)_\infty}.$$
Finally, observing by definition that $\dtli = (q; \, q)_\infty^l \, (tq;
\,q)_\infty^{-l}$ and putting everything into equation \eqref{puttog},
theorem \ref{thm-lev1-tsf} is proved. \qed

Equation \eqref{strid} immediately gives:
\begin{corollary}
 Let $\kma$ be one of the simply laced
  untwisted affine Lie algebras $A_l^{(1)}, D_l^{(1)}, E_l^{(1)}$. Then
$$\alolo =  \frac{1}{\dprod_{i=1}^l (t^{d_i} q;\,q)_{\infty}} =
    \frac{1}{\dprod_{i=1}^l \dprod_{n=1}^\infty (1-t^{d_i} q^n)}$$
where $d_1, d_2, \cdots, d_l$ are the degrees of the underlying finite
    dimensional simple Lie algebra $\kmaf$ ($=A_l, D_l, E_l$).
\end{corollary}

Comparing coefficients of $q=e^{-\delta}$ in \eqref{thmeq} and using
 $K_{\Lambda_0,\Lambda_0- \delta}(1)=l$,  we get another corollary:
\begin{corollary}\label{degcor}
 Let $\kma$ be an untwisted affine Lie algebra. Then,
$$ K_{\Lambda_0,\Lambda_0- \delta}(t) = \sum_{i=1}^l t^{d_i},$$
where $d_i$ are the degrees of $\kmaf$.
\end{corollary}

A strengthened version of this corollary will be proved in \S \ref{positi}.
We remark that this corollary is analogous to the classical result of
Hesselink \cite{hess} and (independently) Peterson \cite{peterson}
which states that $\Kf_{\highroot,0}(t) = \displaystyle\sum_{i=1}^l
t^{e_i}$  where $\Kf$ denotes a Kostka-Foulkes polynomial for $\kmaf$ and
$\highroot$ is the highest long root of $\kmaf$.

\section{On positivity of $\kost$}\label{positi}
For finite dimensional simple Lie algebras, the classical
theory  shows that the  Kostka-Foulkes polynomials are essentially certain
Kazhdan-Lusztig polynomials associated to the corresponding extended
affine Weyl group. Via their geometric interpretation \cite{kl}, one obtains the
non-negativity of their coefficients.

Let $\kma$ be a fixed untwisted affine Kac-Moody algebra with
underlying finite dimensional simple Lie algebra $\kmaf$.
In this section, we will relate the Kostka-Foulkes polynomials $\kos$
of $\kma$ to classical Kostka-Foulkes polynomials
$\Kf_{\beta,\gamma}(t)$ associated to $\kmaf$.

If $V$ is a $\kmaf$-module and $\gamma$ is in the weight lattice $P(\kmaf)$ of
$\kmaf$, define $\Kf_{V,  \gamma}(t)$  in the natural manner: decompose $V = \oplus_\pi
(\Lf(\pi))^{\oplus m_\pi}$ into a direct sum of irreducible $\kmaf$ modules and set 
$\Kf_{V, \gamma}(t) := \sum_\pi m_\pi \,\Kf_{\pi, \gamma}(t)$. Given a weight
$\lambda \in P(\kma)$, let $\lf$ denote its restriction
to the Cartan subalgebra of $\kmaf$; clearly  $\lf \in P(\kmaf)$.

The main result of this section is the following.
\begin{proposition}\label{post}
Let $\kma$ be an untwisted affine Kac-Moody algebra and let $\lambda
\in P^+(\kma)$ such that $(\lambda, \alpha^\vee_0) \geq 1$. Then
$K_{\lambda, \lambda -\delta}(t) = t\, \Kf_{V, \lf}(t)$ where $V:=
\kmaf \otimes \Lf(\lf)$.
\end{proposition}
Before giving the proof, we mention the following immediate corollaries:
\begin{corollary}
For $\kma, \lambda$ as in the proposition, we have 
$K_{\lambda, \lambda -\delta}(t) \in \nat[t]$.
\end{corollary}
We note that  it is not  true that $\kost \in
\nat[t]$ for {\em all pairs} of dominant weights $\lambda, \mu$. As an
example, let $\kma = A_1^{(1)}$ (affine $sl_2$), $\lambda =0$, $\mu =
-\delta$; then an easy calculation using theorem \ref{mainthm} gives
$\kost = t^2 - t$. 
\begin{corollary}
Let $\kma$ be an untwisted affine Kac-Moody algebra (not necessarily
simply laced). Then for each $p
\geq 1$,
$K_{p \Lambda_0,p \Lambda_0- \delta}(t) = \sum_{i=1}^l t^{d_i}$ where the
$d_i$ are the degrees of $\kmaf$.
\end{corollary}
This follows since $\lambda = p\Lambda_0$ implies that  $\bar{\lambda}=0$
and $ V= \kmaf$. We also note this is a strengthening of corollary \ref{degcor}.
We will now prove proposition \ref{post}. The proof is a
straightforward (if lengthy) calculation.

\noindent
{\em Proof of proposition \ref{post}}:
We assume that the underlying finite dimensional simple Lie algebra
$\kmaf$ has rank $l$, highest long root $\highroot$ and Weyl group $\wf$.
 We have
$$K_{\lambda, \lambda -\delta}(t) = [e^0]  \quad \frac{\sum_{w \in W} (-1)^{\len(w)} e^{w
  (\lambda + \rho) - (\lambda  + \rho)}\cdot e^\delta}{\prod_{\alpha
  \in \rp} (1-te^{-\alpha})^{m_\alpha}}.$$

If $w \in W$ has nonzero contribution to the sum on the
  right, then it must satisfy $(\lambda+\rho) - w(\lambda+\rho) \leq
  \delta$.

\vspace{2mm}
\noindent
{\em Claim:} This implies that  $w \in \wf$.\\
{\em Proof:} Write $\Lambda := \lambda + \rho$ and let $w = r_{i_1}
  r_{i_2} \cdots r_{i_k}$ be a reduced word for $w$. Here, $r_{i_j}$
  are simple reflections in $W$ ($0 \leq i_j \leq
  l$). Observe that 
\beq \label{lwl}
\Lambda - w\Lambda = \sum_{p=1}^k (\Lambda,  \alpha^\vee_{i_p}) \beta_p ,
\eeq
where $\beta_p := r_{i_1} r_{i_2} \cdots r_{i_{p-1}}(\alpha_{i_p}) \in
  \rp$ for all $j$. Suppose $ w \not\in \wf$, let $j$ be the least
  index such that $i_j=0$. Then $\beta_j = \alpha_0 + \alf$ where
  $\alf$ is a non-negative integer linear combination of $\alpha_{i_1},
  \cdots, \alpha_{i_{j-1}}$. Now $(\Lambda,  \alpha^\vee_{i_j}) =
  (\lambda + \rho,  \alpha^\vee_{0}) \geq 2$. Thus the $j^{th}$ term
  on the right hand side of equation \eqref{lwl} is $\geq \,2(\alpha_0 +
  \alf)$. This is clearly $\nleq \delta$, since the coefficient of
  $\alpha_0$ in $\delta$ is 1. This proves our claim. \qed

Further, $w \in \wf$
  implies that  $w\lambda - \lambda = w \lf - \lf$. This implies
 that  $K_{\lambda, \lambda -\delta}(t)$ is the coefficient of  $e^0$ in
  the product
$$\frac{\dsum_{w \in \wf} (-1)^{\len(w)} e^{w
  (\lf + \rho) - (\lf  + \rho)}}{\dprod_{\beta
  \in \delp} (1-e^{-\beta})} 
\; \prod_{\beta
  \in \delp} \frac{1-e^{-\beta}}{1-t e^{-\beta}} 
\; \; \frac{e^\delta}{\dprod_{\alpha
  \in \rp \backslash \delp} (1-te^{-\alpha})^{m_\alpha}}.$$
By the Weyl character formula, the first term equals
  $e^{-\lf}\schf[\lf]$. Put
\begin{align}
\xi &:= \prod_{\beta \in \delp}
\frac{1-e^{-\beta}}{1-te^{-\beta}} = \prod_{\beta \in \delp}(1 +
(t-1)e^{-\beta} + (t^2-t)e^{-2\beta} + \cdots) \text{ and } \label{xiq}\\
\eta &:= \frac{e^\delta}{\displaystyle\prod_{\alpha
    \in \rp \backslash \delp} (1-te^{-\alpha})^{m_\alpha}}
=  e^\delta \displaystyle\prod_{\alpha\in \rp \backslash \delp} (1 +
  m_\alpha t \;e^{-\alpha} + \binom{m_\alpha+1}{2} t^2 \; e^{-2\alpha} +
  \cdots). \label{etaq}
\end{align}

 Thus
$K_{\lambda, \lambda -\delta}(t) = [e^0]  \; (e^{-\lf}\schf[\lf]) \;\xi
  \eta$. Note that $ (e^{-\lf}\,\schf[\lf]) \,\xi$ is a power series only
  involving $e^{-\gamma}, \gamma \in \Qf^+$. We analyze $\eta$ more closely. Given $\alpha \in Q$, write $\alpha =
  \sum_{i=0}^l n_i(\alpha) \alpha_i$ with $n_i(\alpha) \in \integers$. Then
  $\alpha \in \rp \backslash \delp$ implies that  $n_0(\alpha) \geq 1$.
Consider  the expression on the far right of equation \eqref{etaq}.
For each $\alpha\in \rp \backslash \delp$, 
suppose we throw away the terms $e^{-p\alpha} \; (p \geq 1)$ for which
$p\alpha \nleq \delta$, it is clear that this does not affect the
value of $K_{\lambda, \lambda -\delta}(t)$.
If $n_0(\alpha) \geq 2$, we can throw away $e^{-p\alpha}$ for all $p \geq 1$ while if $n_0(\alpha) =1 $, we can throw away $e^{-p\alpha}$ 
for all $p \geq 2$. This means that if we define
$$\tilde{\eta} := e^\delta \prod_{\substack{\alpha\in \rp \\ n_0(\alpha)=1}} (1 +
  m_\alpha t \;e^{-\alpha}),$$
then the coefficients of $e^0$ in $ (e^{-\lf}\schf[\lf]) \xi \eta$ and
  $(e^{-\lf}\schf[\lf]) \xi \tilde{\eta}$
  are the same. Observe now that $\{\alpha\in \rp : n_0(\alpha)=1\} =
  \{\delta\} \cup \{-\beta+\delta: \beta \in \delp\}$. Since $m_\delta
  =l$ and $m_{-\beta + \delta} = 1$, this implies that 
$$\tilde{\eta} := e^\delta (1+lte^{-\delta}) \prod_{\beta\in \delp} (1
+  t e^{\beta-\delta}).$$
Expanding out, we get
$$ \tilde{\eta} := lt + t\sum_{\beta \in \delp}
  e^{\beta} + \sum_{\substack{\gamma \in Q\\ n_0(\gamma) \neq
  0}} p_\gamma (t) e^\gamma$$
for some $p_\gamma(t) \in \integers[t]$, i.e, the last sum runs over
  $\gamma \in Q\backslash \Qf$. As remarked before,
  $(e^{-\lf}\schf[\lf]) \xi$ is a power series involving
  $e^{-\gamma}$ with $n_0(\gamma)=0$. Thus, finally
\beq \label{ptone}
K_{\Lambda_0,\Lambda_0- \delta}(t) = lt + t\sum_{\beta \in \delp}
  [e^{-\beta}] \: (e^{-\lf}\schf[\lf])\, \xi\,.
\eeq

We now turn to $\kmaf$. We have  $V=\Lf(\highroot)
\otimes \Lf(\lf)$. This decomposes into
$$\Lf(\highroot) \otimes \Lf(\lf) = \bigoplus_{\pi \in \Pf^+} c_\pi
\Lf(\pi) \;\; \text{ with } c_\pi \in \integers^{\geq 0}.$$
Thus $\schf[\highroot] \schf[\lf] = \sum_\pi c_\pi
\schf[\pi]$. Consider the sum 
\begin{align}
\sum_{\pi \in \Pf^+} c_\pi \Kf_{\pi,\lf}(t) &= [e^0] \quad \frac{\sum_\pi c_\pi \sum_{w \in
    \wf} (-1)^{\len(w)} e^{w (\pi + \rho) -  (\lf+\rho)}}
{\prod_{\alpha \in \delp} (1-te^{-\alpha})} \notag\\
&= [e^0]\quad e^{-\lf} \: \xi \: \sum_\pi c_\pi \schf[\pi] \notag\\
&=[e^0]   \quad e^{-\lf} \,\xi \,\schf[\highroot] \,\schf[\lf]\notag\\
&= [e^0] \quad (e^{-\lf}\schf[\lf] \:\xi) (l+ \sum_{\alpha \in \delp}
  (e^{\alpha} + e^{-\alpha})) \notag\\
&= l + \sum_{\alpha \in \delp}[e^{-\alpha}] \; (e^{-\lf}\,\schf[\lf]\, \xi) \label{secpt}.
\end{align}
Comparing equations \eqref{ptone} and \eqref{secpt}, proposition
\ref{post} is proved. \qed

\vspace{1.5mm}
\noindent
{\bf Remark:} The above proof works with no essential change if we
replace the node 0 by a node $p$ which satisfies the property that if
$\delta = \sum_{j=0}^l a_j \alpha_j$, then $a_p=1$.
Since all nodes of the Dynkin diagram of
$A_n^{(1)}$ satisfy this property, we have the following corollary.
\begin{corollary}
Let $\kma = A_l^{(1)}$, and $\lambda$ be a dominant weight such that not all
$(\lambda, \alpha^\vee_i)$ are zero (equivalently $\lambda \not\in \complex
\delta$). Then $K_{\lambda, \lambda - \delta}(t) \in \nat[t]$.
\end{corollary}

\section{Miscellaneous facts}\label{misc}
In this subsection, we collect together an assortment of facts
concerning the Hall-Littlewood functions $\hl$. We work in the full
generality of $\kma$ being a symmetrizable Kac-Moody algebra.

\vspace{2mm}
\noindent
7.1. Recall from \eqref{plclchim} that $\hl = \sum_{\mu \in P^+}
   c_{\lambda\mu}(t) \: \sch[\mu]$. The $c_{\lambda\mu}(t)$ are the
   entries of the inverse matrix of $[\kost]_{\lambda,\mu}$. Since
   this latter matrix is lower triangular with ones on the diagonal,
   each $c_{\lambda\mu}(t)$ is a polynomial in the $K_{\pi\gamma}(t)$
   with integer coefficients. Theorem \ref{mainthm} then implies that
   $c_{\lambda\mu}(t) \in \integers[t]$ too (note that \eqref{corr9}
   only guaranteed $c_{\lambda\mu}(t) \in
   \integers[[t]]$). This means we can specialize $t$ to any complex
   number $a$ and $P_\lambda(a)$ would end up being a well-defined
   element of $\Ep$.

\vspace{1.5mm}
\noindent
7.2. When $\kma$ is finite dimensional, it follows from an argument due
   to Stembridge \cite{jrs} 
   that $c_{\lambda\mu}(t)$ has the following expression: 
\beq \label{clmexp}
c_{\lambda\mu}(t) := \underset{w(\lambda+\rho - |A|) = \mu + \rho}{\sum_A \sum_{w \in W}}
      (-1)^{\len(w)} (-t)^{\#A},
\eeq
where the outer sum only ranges over those $A \subset \rp$ for which
 $(\lambda,\alpha) >0$ for all $\alpha \in A$ ({\em cf} equation
   \eqref{corr9}). This expression
   explicitly demonstrates that $c_{\lambda\mu}(t)$ is a polynomial in
   $t$. Stembridge's elegant argument can be adapted     without
   change to our present situation; in our case, it shows that for any
 \skma $\kma$ and $\lambda, \mu \in P^+$ such that $\# W_\lambda <
   \infty$, $c_{\lambda\mu}(t)$ is given by the same (explicitly
   polynomial) expression of equation \eqref{clmexp}, where we now let
   the outer sum range over all multisets $A \in \scA$ which satisfy
   $(\lambda, \alpha) > 0$ for all $\alpha \in A$.

\vspace{1.5mm}
\noindent
7.3. Next, we prove a necessary condition on $(\lambda,\mu)$
   in order that $c_{\lambda\mu}(t) \neq 0$.

\begin{proposition}
Let $\kma$ be a \skma, $\lambda, \mu \in P^+$. Then $c_{\lambda\mu}(t)
\neq 0$ implies that there exists $B \in \scA$ such that $\mu = \lambda -
|B|$.
\end{proposition}
\pfbegin
From equation \eqref{corr9}, we have $c_{\lambda\mu}(t) \neq 0$ implies that there exists $w \in W, A \in \scA$ such that $\mu + \rho = w(\lambda) + w(\rho
- |A|)$. Since the formal character of $L(\rho)$ is
$e^\rho\prod_{\alpha \in \rp} (1 + e^{-\alpha})^{m_\alpha}$, the weights
of $L(\rho)$ are $\{ \rho - |B|: B \in \scA\}$. Let $\gamma \in P^+$
be the unique dominant weight in the Weyl group orbit of $\rho -
|A|$. Then $\gamma \leq \rho$ and there exists $\tau \in W$ such that
$w(\rho  - |A|) = \tau \gamma$. So, we have $\mu + \rho = w\lambda +
\tau \gamma$ with $\lambda, \gamma, \mu+\rho \in P^+$. By the
Parthasarathy-Ranga Rao-Varadarajan (PRV) conjecture (see \cite{L1}
for a recent proof), this implies that $L(\mu+\rho)$ occurs as a
summand in the decomposition of $L(\lambda) \otimes L(\gamma)$. It is
well-known (see for eg \cite{L2}) that this implies $\mu +
\rho = \lambda + \eta$ for  a weight $\eta$ of $L(\gamma)$.\\
{\em Claim:} $\eta$ is also a weight of $L(\rho)$\\
{\em Proof:} Let $\eta^+$ be the dominant weight in the Weyl orbit of
$\eta$. Then $\eta^+ \leq \gamma \leq \rho$. 
Since $(\rho,\alpha_i) >0$ for all $i$, proposition {11.2}
 of \cite{kac} implies that $\eta^+$ is
in fact a weight of $L(\rho)$. Hence, so is $\eta$. 

This of course means that there exists $B \in \scA$ such that $\eta = \rho
- |B|$. Thus $\mu + \rho = \lambda + \rho - |B|$ proving our proposition.
\qed

We now deduce a corollary that was proved by R. Brylinksi in \cite{rkg}.
\begin{corollary}
Let $\kma$ be finite dimensional. Then $c_{\lambda\mu}(t) \neq 0$
implies that $\mu \geq \lambda - 2 \rho$.
\end{corollary}
\pfbegin
By our proposition above, $\mu = \lambda - |B|$ where $B$ is a subset
of $\rp$. This means $|B| \leq \sum_{\alpha \in \rp} \alpha = 2 \rho$. \qed

\vspace{2mm}
\noindent
7.4. Let $\kma$ be a \skma and consider the special case that $\lambda
   \in P^{++}$ (a regular dominant weight) i.e, $(\lambda, \alpha_i)
   >0$ for all $i=1\cdots n$. Then $W_\lambda(t)=1$. Consider the
   specialization of $\hl$ at $t=-1$ (for classical type A, 
these are the Schur Q-functions). Now,
$$P_\lambda(-1) = \frac{\sum_{w \in W} (-1)^{\len(w)} w(e^{\lambda+\rho}
\prod_{\alpha \in \rp} (1 + e^{-\alpha})^{m_\alpha})}{e^\rho
\prod_{\alpha \in \rp} (1 - e^{-\alpha})^{m_\alpha}}.$$
As mentioned before, $\sch[\rho] = e^\rho \prod_{\alpha \in \rp} (1 +
e^{-\alpha})^{m_\alpha}$ is $\ch L(\rho)$ and thus $W$-invariant. Thus
$$ P_\lambda(-1) = \frac{J(e^\lambda)}{J(e^\rho)} \sch[\rho].$$
By the Weyl-Kac character formula, this is just $\sch[\lambda-\rho]\,
\sch[\rho]$ (observe that $\lambda   \in P^{++}$ ensures $\lambda-\rho
\in P^+$). Since $c_{\lambda\mu}(-1)$ is the coefficient of
$\sch[\mu]$ in $P_\lambda(-1)$, we obtain the following generalization of a well known result in the classical case.
\begin{proposition}
$\lambda \in P^{++}, \mu \in P^+$ implies that  $c_{\lambda\mu}(-1)$ is the
  multiplicity of $L(\mu)$ in the tensor product
  $L(\lambda-\rho)\otimes L(\rho)$. Thus $c_{\lambda\mu}(-1) \in
  \integers^{\geq 0}$
\end{proposition}

\end{document}